\documentclass[12pt]{amsart}
\usepackage{amsmath,amsfonts,amssymb}

\def\le{\leqslant}

\def\ge{\geqslant}

\def\frak{\mathfrak}

\def\cH{\mathcal{H}}

\def\cM{\mathcal{M}}

\numberwithin{equation}{section}

\newtheorem{cor}{COROLLARY}[section]
\newtheorem{prop}[cor]{PROPOSITION}

\theoremstyle{definition}

\begin{document}
\title{ When Euler met Brun}

\author[Friedlander]{John B. Friedlander}

\maketitle

\dedicatory \quad\quad\quad\quad\quad\, \,{\sl For the 75th birthday of Henryk Iwaniec}

\medskip

{\bf Abstract:}
 In this note we survey and study some aspects of the distribution of primes in very short intervals. 
\footnote{key words: primes, sieves, reciprocals, short intervals}

\section{\bf Introduction}

L. Euler famously proved that the sum of the reciprocals of the prime numbers is divergent, strengthening a millennia old theorem of Euclid. It is reasonable to view Euler's result, at least when taken together with Dirichlet's profound generalization to primes in arithmetic progressions, as having launched the modern subject of analytic number theory. 

V. Brun famously proved that, when restricted to a sum over the reciprocals of twin primes, the series is in contrast convergent, along the way launching the modern subject of sieve theory. 

If we let $p'$ denote the least prime exceeding the prime $p$, its successor prime, then Brun's theorem may be stated as $\sum_{p' =p  +2}1/p < \infty$. 
It is well-known and clear from the proof that, with virtually no modification, Brun's theorem holds as well for the sum $\sum_{p' = p+ k}1/p$ for any even integer $k$ and, since it is trivial for odd $k$, we have
\begin{equation*}
\sum_{p'-p \le K}\frac{1}{p} < \infty ,
\end{equation*}
for any fixed $K$. 

So a question that  naturally suggests itself is to wonder what happens when we replace $K$ by a growing function, say $y(p)=\lambda (p) \log p$, and try to study how quickly growing a function we can take and still maintain the conclusion
\begin{equation}\label{eq:1.1}
\lim_{x \rightarrow \infty}\sum_{\substack{p\le x\\p'-p \le y(p)}}\frac{1}{p} < \infty .
\end{equation}

We know, by a result of Goldston, Pintz and Yildirim, Theorem 1 of [GPY4], that for any fixed $\lambda > 0$, a positive proportion $c(\lambda)$ of the primes $p\le x$ have successor satisfying the bound $p'- p < \lambda \log p$. Thus we are not able to take $y(p) = \lambda\log p$ for any fixed $\lambda > 0$, no matter how small. This means that we are dealing with prime gaps that are of smaller order of magnitude than the average. 

On the other hand, if we could take a choice of $y$ which yielded, for some fixed real $a>1$, an upper bound 
\begin{equation*}
\sum_{\substack{p\le x\\p'-p \le y(p)}}1 \ll \frac{x}{(\log x)^a}
\end{equation*}
then partial summation would show the convergence of its sum of reciprocals. Recall that, in the case of twin primes, we have such a bound with the comfortable margin $a=2$. Moreover, for the full set of primes ($y=p$) the divergence, namely $\sum_{p\le x} 1/p \sim \log\log x$, just barely holds. This suggests that the set of twin primes is substantially further from the ``meeting place'' and that the relevant gaps are short but not so terribly short. 

After a tiny bit of history and a few remarks in Section 2, we introduce in Section 3 what we believe to be a reasonable conjecture which allows us to deduce fairly close upper and lower bounds for the function $y(p)$ in question. Then, in Section 4, we consider what we are able to learn about the problem unconditionally using sieve bounds. The proofs are elementary. 
 
{\vskip 0.3 in}

{\bf Acknowledgements:} Thanks go to the referee for pointing out a mistake in an earlier version. Research of the author supported in part by an NSERC Canada Research Grant, continuously since 1981. The problem considered in this paper was inspired by a question from Peter Rosenthal, colleague since 1980. Title suggested by a work of Robert Reiner. This paper was written in honour of Henryk Iwaniec, frequent collaborator and constant friend since 1976.

\section{\bf Small gaps and very small gaps} 

By the Prime Number Theorem we know that the average gap between primes of size around $x$ is $\log x$. As it happens however, our knowledge of primes is still sufficiently limited that, for at least one problem, that of maximal gaps, intervals as large as $x^{.49}$ are still too small to be within our reach. Nevertheless, for purposes of this paper, we think of prime gaps near $x$ having length $\lambda\log x, \lambda <1$ a constant, as being ``smaller than average'', those  unbounded but of size $y(x)=o(\log x)$ as being ``very small'' and of bounded gaps as being, well, bounded. 

At the two extremes of these three gap ranges there has, after many years, been phenomenal recent progress. Thus, for the smaller than average gaps $\lambda \log x$ one has the breakthroughs of Goldston, Pintz and Yildirim [GPY1], wherein they show that $\lim\inf (p'-p)/\log p =0$ and subsequently as well the positive proportion result already mentioned. 

The GPY bound was greatly strengthened to encompass the bounded gaps by the work of Zhang [Zha], then Maynard [M1], who showed that $\lim\inf (p'-p)$ is bounded and, in the latter case, quite a bit more as well. 

The functions $y(x)$ which are near the meeting place for our problem lie in what we have called very small gaps, the largely uncharted region between these two extreme interval lengths. 

There seems to be not a great deal known about this middle range. One exception comes from Maynard's sieve which, although originally designed to treat distribution in bounded gaps, still has striking consequences which venture into the middle range. Thus, Theorem 3.2 of Maynard [M2] states:

{ \vskip .1 in} 

{\bf Theorem (Maynard):}
For any $x, y \ge 1$ there are $\gg x\exp (-\sqrt{ \log x})$ integers $x_0\in [x,2x]$ such that 
\begin{equation*}
\pi(x_0 +y) -\pi(x_0) \gg \log y .
\end{equation*}

{\vskip .1 in }

Another exception, closer to our problem, is the sieve-driven upper bound, Theorem 2 of [GPY4], the statement of which we postpone until Section 4.

Finally, we mention the papers [EN] of Erd\"os and Nathanson and of Zhou [Zho] which study a problem of flavour similar to ours but dealing with a sum over (suitably weighted) reciprocals of prime gaps.

\section{\bf Conditional results} 

Despite the breakthroughs of GPY and Maynard, much of what is thought to be true about the distribution of small prime gaps is either conjecture or else deductions from conjectures. Thus, Granville and Lumley [GrLu] study the maximum number of primes in a gap near $x$ having length $y$ and conjecture that it is 
asymptotic to $y/\log y$ ``for $y\le c \log x$ as long as $y\rightarrow \infty$ as $x\rightarrow \infty$''. 

Nevertheless, most of what appears in the literature on gaps of below average length seems to deal only with gaps of size $ \lambda \log x$ with constant $\lambda$ rather than $\lambda =o(1)$. A pioneering work which is concerned with this range of $\lambda$ constant is that of Gallagher [G]. 

Gallagher's theorem on primes in short intervals rests on the assumption of a uniform version of a well-known conjecture due to Hardy and Littlewood [HL].

Let $\cH = h_1,\ldots , h_r$ denote a set of distinct non-negative integers and let $\nu_{\cH} (p)$ denote the number of residue classes modulo $p$ occupied by the $h_j$. Assume that the set $\cH$ is ``admissible'' in the sense that $\nu_{\cH} (p) <p$ for all primes $p$. We recall the definition of the relevant``singular series''  

\begin{equation}\label{eq:3.1}
{\frak S(\cH)} =\prod_p\left(1-\frac{\nu_{\cH} (p)}{p}\right)\left(1-\frac{1}{p}\right)^{-r} .
\end{equation}

{\bf Hypothesis G:}
Let $\pi(x;\cH)$ denote the number of integers $n\le x$ such that $n+h_j$ is prime for all $h_j \in \cH$. Then 
\begin{equation}\label{eq:3.2}
 \pi(x;\cH) \sim {\frak S(\cH)}\frac{ x}{(\log x)^r}, 
\end{equation}
as $x\rightarrow\infty$ uniformly for all $h_1,\ldots , h_r \le \lambda \log x$, where $\lambda >0$ is fixed. 

A similar conjecture (on the number of prime couples and prime triples with fixed spacing) was used by Goldston and Ledoan [GoLe2] in relation to their work on the prime jumping champions.

{\vskip .1 in}

{\bf Theorem (Gallagher):}
  Denote  by  $P_k(h,x)$  the  number  of  integers  $n \le x$ for which the
interval  $(n, n + h]$  contains exactly $k$ primes. Then, under the above Hypothesis G, with $h= \lambda \log x$, $\lambda> 0$ fixed, 
\begin{equation}\label{eq:3.3}
 P_k(h,x)\sim xe^{-\lambda}\frac{\lambda^k}{ k!}
\end{equation}
as $x\rightarrow\infty$.  

{\vskip .1 in}

Soundararajan pointed out, in very slightly different form [S, exercise 1.3] (see also Goldston and Ledoan [GoLe1]), the following consequence of Gallagher's theorem.

{\vskip .1 in}

{\bf Corollary:} Again under the Hypothesis G of the theorem, we have for fixed 
$\lambda >0$, 
\begin{equation*}
\frac{1}{\pi(x)}\#\left\{p\le x; \frac{p'-p}{\log p} \le \lambda\right\}
\sim \int^{\lambda}_{0}e^{-t}dt = 1-e^{-\lambda}.
\end{equation*}
 
{\vskip .1 in}

In order to close in on a consequence of this for our problem, we shall assume that this statement holds in a range slightly different from $\lambda$ fixed, but rather with $\lambda \rightarrow 0$ slowly as $x \rightarrow \infty$. This seems reasonable, especially if the approach is sufficiently slow, for instance the following is more than ample. 

{\vskip .1 in}

{\bf Hypothesis S:} Assume that, as $t \rightarrow \infty$ with $\lambda (t) \rightarrow 0$, subject to 
\begin{equation}\label{eq:3.4}
\lambda(t) \gg \frac{1}{(\log\log t)^2}  
\end{equation}
we have
\begin{equation}\label{eq:3.5}
\frac{1}{\pi(x)}\#\left\{p\le x; \frac{p'-p}{\log p} \le \lambda(p)\right\}
\sim \int^{\lambda(x)}_{0}e^{-u}du = 1-e^{-\lambda(x)}.
\end{equation}

{\vskip .1 in}

\noindent
Note that, for $\lambda = o(1)$ we have $1- e^{-\lambda} \sim \lambda $.

We prefer to base our conditional statements on Hypothesis S which, as it concerns a {\it sum} over gaps of specific lengths, one might believe it could hold even if the previous two conjectures, which deal with individual spacings, did not. We are not suggesting that Hypothesis S holds without further conditions on the function $\lambda$; we only require it for the specific functions mentioned in Proposition 3.1 below. 

As it happens, we shall need iterated logarithms to enter our statements so we define, as usual, for $k\ge 2$, $\log_k$ to be the $k$-th iterated logarithm: 
\begin{equation}\label{eq:3.6} 
\log_k = \log\log_{k-1} , \quad \log_1 = \log 
\end{equation}
and what we might call the logorial function
\begin{equation}\label{eq:3.7} 
Log_k = \prod_{2\le j \le k}\log_j .
\end{equation}
For each of these functions, we can ignore the finite set of small integers where it is not defined. 

Applying partial summation to the sum in~\eqref{eq:1.1}, then Hypothesis S, next inputting in the resulting integral the Prime Number Theorem, then making the change of variable $u= \log_{k+1}t$, respectively $u= \log_kt$, one deduces the following results which narrow in on the answer to our question. 

\begin{prop}
Let $k\ge 2$ be a fixed integer and $\varepsilon >0$ a fixed real. Denote  
$y(p)=\lambda(p) \log p$. Then we have as $x \rightarrow \infty$
\begin{equation}\label{eq:3.8}
\sum_{\substack{p\le x\\p'-p \le y(p)}}\frac{1}{p} \sim \log_{k+1}x \rightarrow\infty \quad { for }
\,\, \lambda (p) = 1/ Log_k(p) ,
\end{equation}
but
\begin{equation}\label{eq:3.9}
\sum_{\substack{p\le x\\p'-p \le y(p)}}\frac{1}{p} \quad  { is\,\, bounded\,\, for} 
\,\, \lambda (p) = 1/ Log_k(p)(\log_k p)^{\varepsilon},
\end{equation}
under the assumption that~\eqref{eq:3.5} holds for these functions $\lambda$.
\end{prop} 

In fact, as we shall see in the next section, the latter statement~\eqref{eq:3.9} holds unconditionally, although the former seems hopeless with current tools.
We may remark that Theorem 1 of [GPY4] gives a lower bound for the left hand side of~\eqref{eq:3.5}, weaker but unconditional and holding for intervals of smaller than average size, that is $\lambda$ constant. If this could be extended to also hold for intervals having length in a suitable part of the range $o(\log x)$, perhaps using ideas from [M2], this might give a correpondingly weaker but unconditional result in the direction of~\eqref{eq:3.8}. 

Note that~\eqref{eq:3.8}, if true, is not the end of the story since we could then construct more artificial choices for $y(p)$ which are closer to the breakpoint. To see this, begin with $\lambda = 1/Log_k$ for some particular $k$ as in~\eqref{eq:3.8}, carry on until the sum of reciprocals of primes up to $x$ exceeds $1$ (still more complicated choices go further) and also $x$ is large enough that $\log_{k+1}x$ is defined. Then, for ensuing $p$, replace $1/Log_k$ by $1/Log_{k+1}$ in the definition of $y(p)$, continue summing  and iterate.

\section{\bf Sieve survivors}

By a sieve ``survivor'' we mean what is more awkwardly called ``an integer with no small prime factor'' but in recent years has, at least in sieve theory, become confused with an ``almost-prime'' even though the latter is invariably defined as an integer with few prime factors. There are many sieve results which produce what are called almost-primes when they actually produce integers from the special subset of survivors. When we use the term here we are thinking of an integer $\le x$ having no prime factor $<z$ with $z =x^{\delta}$ for some $\delta >0$, possibly small but fixed. The survivors, at least in one sense, more strongly resemble the primes, occurring as they do with the same order of magnitude. Recall that the almost-primes occur with a greater order of magnitude, how much so depending on the number of prime factors permitted.

In proving lower and upper bounds for the sum $\sum_p1/p$ by partial summation, it is obvious that we do not need an asymptotic formula such as that conjectured in~\eqref{eq:3.5}, but only sufficiently good lower, respectively upper, bounds for the corresponding prime counting function.  This suggests that sieve methods can prove useful. Already, in the work of Gallagher [G], his Theorem 2 gives an upper bound in connection with his question.

More recently, again in the case of the upper bound, one has a result more closely related to our topic, namely Theorem 2 of Goldston, Pintz and Yildirim [GPY4] which states: 

For any $h>2$ as $x\rightarrow\infty$ we have
\begin{equation}\label{eq:4.1}
\sum_{\substack{p\le x\\p'-p \le h}} 1 \ll \min\{h/\log x, 1\} \pi (x),
\end{equation}
which for us is of interest when $h=o(\log x)$, although it is also meaningful if $h \le c\log x$ provided that the constant $c$ is sufficiently small. This shows that, if in~\eqref{eq:3.5} we replace the conjectured asymptotic by an upper bound, then that result is unconditionally true in a wide range.

{\vskip 0.1 in}
We shall take $z=x^{\delta}$ for some positive $\delta$. Let $m$ run through $\cM$, the set of integers $m\le x$ which are the survivors of sieving by the primes less than $z$. For each $m$ let $m'$ be its successor, the first larger integer in the set.
\vskip .1 in

A proof, perhaps slightly streamlined, of~\eqref{eq:4.1} can proceed as follows. We start with a sieve upper bound (e.g. Theorem 7.16 of [Opera]):
\begin{equation}\label{eq:4.2}
\sum_{\substack{m\le x\\m, m+d \in \cM}} 1 \ll {\frak S}_d  \frac{x}{(\log x)^2} , 
\end{equation} 
which holds for every integer $d,\, 1\le d\le h$ with the same implied constant depending only on $\delta$. Here, ${\frak S}_d $ is the singular series~\eqref{eq:3.1} for the set $\cH =\{0, d\}$. 
Summing over $1\le d\le h$ and using the fact that 
\begin{equation}\label{eq:4.3}
\sum_{d\le h}{\frak S}_d \sim h 
\end{equation} 
as $h\rightarrow\infty$ 
(much more is known; see [FG], Proposition 1), 
we obtain a bound
\begin{equation}\label{eq:4.4}
\sum_{\substack{m\in\cM\\ m'-m \le h}} 1 \le\sum_{\substack{m_1, m_2\in\cM\\ 1\le m_1-m_2 \le h}} 1 \ll h  \frac{x}{(\log x)^2} \sim \frac{h}{\log x}\pi (x). 
\end{equation} 
This bound applies, a fortiori, to the corresponding sums over primes since the number of primes less than $z$ offers a negligible contribution. This proves~\eqref{eq:4.1}.

\vskip .1in

Note that the upper bound in~\eqref{eq:4.4} is actually a bound for a larger number of $m\in\cM$, specifically those with two {\it or more} unsifted integers in the interval $[m, m+h]$. As it happens, the number of triples in such a short interval is, as we shall see, of smaller order of magnitude and we lose nothing for this result by ignoring them. 

On the other hand, returning to the sum in~\eqref{eq:4.2}, we can instead begin with a sieve lower bound. If  we have chosen for instance $\delta = 1/10$, then $z$ is sufficiently small to enable a positive lower bound for this two-dimensional sieve problem (see e.g. Corollary 6.13 of [Opera]): 
\begin{equation}\label{eq:4.5}
\sum_{\substack{m\le x\\m, m+d \in \cM}} 1 \gg {\frak S}_d  \frac{x}{(\log x)^2} .
\end{equation} 
After summing over $d$ we have, in contrast to~\eqref{eq:4.4},  
\begin{equation}\label{eq:4.6}
\sum_{\substack{m_1, m_2\in\cM\\ 1\le m_1-m_2 \le h}} 1 \gg h  \frac{x}{(\log x)^2} \sim \frac{h}{\log x}\pi (x). 
\end{equation} 
Now however, in order to get a lower bound for the number of consecutive pairs corresponding to that in~\eqref{eq:4.4} we need to subtract out (an upper bound for) the contribution of triples in the short intervals. For this, we need an anologue to~\eqref{eq:4.3}. We take a special case of a more general result of
Odlyzko, Rubinstein and Wolf [ORW] which, in turn, extended a basic result of Gallagher [G]; see also Lemma 2 of [GoLe2].

For $\cH= \{0, d_1, d_2\}$ we have
\begin{equation}\label{eq:4.7}
\sum_{d_1,d_2\le h}{\frak S}(\cH)\sim h^2,  
\end{equation} 
as $h\rightarrow\infty$.
This yields an upper bound (again, Theorem 7.16 of [Opera]) 
\begin{equation}\label{eq:4.8}
\sum_{1\le d_1 < d_2 \le h}\sum_{\substack{m\le x\\m, m+d_1, m+d_2\in \cM}} 1 \ll h^2 \frac{x}{(\log x)^3} 
\end{equation} 
and as we intend to take $h=o(\log x)$, this is small compared  to the lower bound in~\eqref{eq:4.6}. Hence we have 
\begin{equation}\label{eq:4.9}
\sum_{\substack{m\in\cM\\ m'-m \le h}} 1 \gg h  \frac{x}{(\log x)^2} \sim \frac{h}{\log x}\pi (x). 
\end{equation} 

We want to use upper and lower bounds for a sum much as that in~\eqref{eq:4.9} but slightly modified. In the first place, we want to change the sifting range so that, rather than having $z=x^{\delta}$ in the definition of $\cM$ we now have $z=m^{\delta}$.  We also want to be able to allow $h$ to depend on $m$, specifically $h=y(m)=\lambda(m)\log m$. To treat this new sum we split the interval $[0,x]$ into dyadic segments $I=(M,2M]$ with perhaps a shorter interval left over. We request that $y(m)$ be a positive slowly increasing function such as either of the two choices in Proposition 3.1. In particular, we want to make use of the fact that $y(m)$ is very nearly constant over every dyadic interval. 

We follow the arguments that led to~\eqref{eq:4.4},~\eqref{eq:4.9}, obtaining the new bounds 
 \begin{equation}\label{eq:4.10}
\sum_{\substack{m\in I\cap\cM \\ m'-m \le y(m)}} 1 \asymp  \frac{y(M)M}{(\log M)^2} .
\end{equation} 
Summing~\eqref{eq:4.10} over the subintervals and again using the fact that $y(m)$ is nearly constant over dyadic intervals, we find the bounds
\begin{equation}\label{eq:4.11}
\frac{1}{|\cM|}\#\left\{m\in \cM; \frac{m'-m}{\log m} \le \lambda(m)\right\}
\asymp  1-e^{-\lambda(x)}.
\end{equation}
Using~\eqref{eq:4.11} in place of~\eqref{eq:3.5} in the argument for Proposition 3.1, we obtain the following result. 

\begin{prop}
Let $k\ge 2$ be a fixed integer and $\varepsilon >0$ a fixed real. Let $\cM$ denote the set of integers $m\le x$ which are free from prime divisors less than $m^{\delta}$ where $\delta >0$ is fixed and sufficiently small that $1/\delta$ exceeds the sifting limit for the two-dimensional (beta or Selberg) sieve. Let $m'$ denote the successor of $m$ in the set $\cM$ and define $y(m)=\lambda(m) \log m$. Then we have, as $x \rightarrow \infty$,
\begin{equation}\label{eq:4.12}
\sum_{\substack{m\le x\\m'-m \le y(m)}}\frac{1}{m} \rightarrow\infty \quad  { for }
\,\, \lambda(m) = 1/ Log_k(m) ,
\end{equation}
but
\begin{equation}\label{eq:4.13}
\sum_{\substack{m\le x\\m'-m \le y(m)}}\frac{1}{m} \quad  { is\,\, bounded\,\, for} 
\,\, \lambda(m) = 1/ Log_k(m)(\log_k m)^{\varepsilon} .
\end{equation}
\end{prop} 

Note that the sieve limit restriction on $\delta$ is not needed in the case of~\eqref{eq:4.13} and that~\eqref{eq:4.13} shows that~\eqref{eq:3.9} is unconditionally true. As far as the lower bound,~\eqref{eq:4.12} lends some additional credence to the belief that~\eqref{eq:3.8} is also true.



\medskip 
\medskip

Department of Mathematics, University of Toronto

Toronto, Ontario M5S 2E4, Canada  \quad (frdlndr@math.toronto.edu)

\end{document}